\documentclass[12pt]{amsart}
\usepackage{amssymb}
\usepackage{latexsym}
\usepackage{amsmath}
\begin{document}

\def\l{{\mathcal L}}
\def\p{{\mathcal P}}
\def\c{{\mathbb C}}
\def\r{{\mathbb R}}
\def\n{{\mathbb N}}
\def\z{{\mathbb Z}}
\def\u{{\mathbb U}}
\def\s{{\mathcal S}}
\def\one{{\mathbb I}}

\newcommand\wt{\widetilde}  
\newcommand\wh{\widehat}
\newcommand\wo{\overline} 

\newcommand{\e}[2]{\begin{pmatrix}{#2}\cr{#1}\end{pmatrix}}

\newtheorem{theo}{Theorem}[section]
\newtheorem{lemm}{Lemma}[section]
\newtheorem{prop}{Proposition}[section]
\newtheorem{defi}{Definition}[section]
\newtheorem{coro}{Corollary}[section]

\title{Quantum automorphism groups of small metric spaces}
\author{Teodor Banica}
\address{Institut de Mathematiques de Jussieu, 175 rue du Chevaleret, 75013 Paris}
\email{banica@math.jussieu.fr}
\maketitle

To any finite metric space $X$ we associate the universal Hopf $\c^*$-algebra $H$ coacting on $X$. We prove that spaces $X$ having at most 7 points fall into one of the following classes: (1) the coaction of $H$ is not transitive; (2) $H$ is the algebra of functions on the automorphism group of $X$; (3) $X$ is a simplex and $H$ corresponds to a Temperley-Lieb algebra; (4) $X$ is a product of simplices and $H$ corresponds to a Fuss-Catalan algebra.

\section*{Introduction}

Let $X$ be a finite set. Consider faithful coactions $v:\c (X)\to\c (X)\otimes H$ of Hopf $\c^*$-algebras $H$ on $X$. If $H$ is commutative, it has to be of the form $\c (G)$, with $G$ a finite group of permutations of $X$. In the general case, a group dual example shows that such an $H$ can be non commutative and infinite dimensional. The universal Hopf $\c^*$-algebra $H(X)$ coacting on $X$ is constructed by Wang in \cite{wa2}.

If $X$ has 1, 2 or 3 points then $H(X)$ is the algebra of functions on the permutation group of $X$. If $X$ has more than 4 points $H(X)$ is non commutative and infinite dimensional. The fusion rules for its irreducible corepresentations are the same as fusion rules for irreducible representations of $SO(3)$, independently of $\# X$.

It is a natural question to ask whether such uniformisation results still hold when adding some extra structure to $X$. The case of coactions on certain fibrations $X\to Y$ is studied in \cite{fc}. A universal Hopf $\c^*$-algebra $H(X\to Y)$ can be constructed and once again, fusion rules for its corepresentations are given by the combinatorics of the Fibonacci graph, independently of the initial data.

At the level of planar algebras, $H(X)$ corresponds to the Temperley-Lieb algebra and $H(X\to Y)$ corresponds to the Fuss-Catalan algebra of Bisch and Jones (\cite{bj1}).

\medskip

In this paper we study the more general situation when we are given a metric $d$ on $X$. A universal Hopf $\c^*$-algebra $H(X,d)$ can be constructed in this case as well. In many cases $H(X,d)$ collapses to the algebra of functions on the usual automorphism group $G(X,d)$. In other cases the coaction of $H(X,d)$ is not transitive, meaning that there exists a non-constant function $f\in\c (X)$ such that $v(f)=f\otimes 1$.

The most interesting case is when $H(X,d)$ is not commutative and its coaction is transitive. For instance in the case where $X$ is a simplex, meaning that all distances between points are equal, $H(X,d)$ is the above $H(X)$. The algebra $H(X\to Y)$ is as well of the form $H(X,d)$, for a suitable choice of the distance $d$.

We prove that for metric spaces $(X,d)$ having at most 7 points the algebras $H(X,d)$ fall into one of these four classes: commutative, non-transitive, corresponding to the Temperley-Lieb algebra, corresponding to the Fuss-Catalan algebra.

\medskip

This is inspired from recent work of Bisch and Jones. In \cite{bj2}, \cite{bj3} they prove the following result. Let $P$ be a spherical $\c^*$-planar algebra satisfying the following conditions. (1) $dim(P_1)=1$ and $P$ is generated by a 2-box $b\in P_2$. (2) $dim(P_2)=3$. (3) $dim(P_3)\leq 13$. Then (4) $P$ is a Temperley-Lieb or a Fuss-Catalan algebra, or is the algebra associated to $\z_3$, or is another algebra, associated to $\z_2$ and $D_5$.

The above result can be reformulated as follows (cf. \cite{wo2}, \cite{pl}). Let $P$ be a spherical $\c^*$-planar algebra satisfying the following conditions. ($1^\prime$) $dim(P_1)=1$ (``transitivity'') and $P$ is generated by a 2-box $b\in P_2$ (``universality''). ($2^\prime$) $P$ embeds into the depth 1 planar algebra associated in \cite{j2} to the inclusion $\c\subset\c (X)$ and $b$ corresponds in this way to a distance matrix (``tannakian conditions''). ($3^\prime$) $\# X\leq 7$. Then ($4^\prime$) $P$ is a Temperley-Lieb or a Fuss-Catalan algebra, or there is a 4-box in $P$ which corresponds to the flip automorphism of $X\times X$ (``commutativity'').

We see that the general conditions (1) and ($1^\prime$) are the same and that the key conditions (2) and ($2^\prime$) are of different nature. The conclusions (4) and ($4^\prime$) look a bit similar just because the auxiliary conditions (3) and ($3^\prime$) are too strong.

The methods we use for the proof are quite different from those of Bisch and Jones. The fact that the 2-box here is given by an explicit operator brings lots of simplifications with respect to the much more general situation considered in \cite{bj2}, \cite{bj3}. As an example, Landau's exchange relation (\cite{la}), which is one of the main tools in classifying planar algebras, says for a 2-box $p$ that the image of the corresponding projection $P\in\l (\c (X))$ must be closed under multiplication, and this is very easy to check.

In fact the proof of the above result just uses our previous work on the Fuss-Catalan algebra \cite{fc} plus some elementary matrix computations. It is also possible to apply the powerful results of Bisch and Jones to the Hopf $\c^*$-algebras $H(X,d)$. After making all necessary identifications, the problem is to check dimension conditions on the planar algebra of $H(X,d)$. The geometric formulation of these conditions requires some combinatorial work, and this is exactly the kind of work we do in this paper.

\medskip

We don't know what happens for $\# X\geq 8$. Some exceptional examples are expected to arise, for instance in connection with complex Hadamard matrices and Krishnan-Sunder biunitary permutation matrices. Associated to such matrices are quite exotic lattice models and commuting squares, which can be described by coactions (cf. \cite{js}, \cite{bh2}, \cite{sq}). In fact, these coactions are the second motivation for the study of the algebras $H(X,d)$. So far the universal Hopf algebras $H(X)$ haven't been useful for getting new results on vertex models, but one can hope that the Hopf algebras $H(X,d)$, which are somehow of ``second level'', could produce better results.

\medskip

From the point of view explained by Bisch and Jones in \cite{bj2} the simplest planar algebras are those generated by a 2-box, and among them, simplest are those having small dimension. There is a similar hierarchy for Hopf $\c^*$-algebras: (I) By Woronowicz's result \cite{wo3} each  unital Hopf $\c^*$-algebra is an inductive limit of pairs $(H,v)$ as in \cite{wo2}. (II) By replacing $v$ with $ad(v\oplus 1)$ we may assume that $v:A\to A\otimes H$ is a coaction on a finite dimensional $\c^*$-algebra. (III) We say that $(H,v)$ has length $k$ if it is isomorphic to the universal pair $(K,w)$ satisfying $Hom(1,v^{\otimes k})\subset Hom(1,w^{\otimes k})$. (This means that the associated planar algebra has planar depth $k$ in the sense of \cite{ls}.) By analogy with discrete groups, each pair $(H,v)$ is a projective limit of length $k$ pairs. (IV) The simplest situation is when $k=2$. Here $(H,v)$ is presented by the relations $dv=vd$, where $d$ is a generator of $Hom(1,v^{\otimes 2})\simeq End(v)$ (the ``2-box''). (V) The simplest case is when $A$ is abelian and $d$ is a metric on its spectrum $X$. Here $H=H(X,d)$.

Generalisations of the algebras $H(X,d)$ can be found by going up through this hierarchy. On the other side, when trying to generalise a metric space, the first thought goes to Connes' notion of a spectral triple (\cite{c}). There are many problems here, the main one being probably to find a good hierarchy of spectral triples, reflecting somehow the corepresentation theory complexity of the associated Hopf $\c^*$-algebra (if any). See Paschke and Sitarz \cite{ps} and Rieffel \cite{r} for more comments on this subject.

\medskip

The paper is organised as follows. 1. is a preliminary section, with material from \cite{wa2} and \cite{fc}. In 2. we establish some basic facts about coactions on finite metric spaces. In 3. we find some geometric rules for detecting commutativity or non-transitivity of $H(X,d)$, the main tool being diagonalisation of distance matrices. In 4. we assume that $\# X\leq 7$ and we prove the main result.

\section{Coactions on finite sets}

Let $X$ be a finite set. We denote by $\c (X)$ the algebra of complex functions on $X$. We use the canonical involution and trace of $\c (X)$, given by the following formulae.
$$f^*(i)=\wo{f(i)}\hskip 2cm tr(f)=\sum_i f(i)$$

If $K$ is a set we denote by $M_X(K)$ the set of functions from $X\times X$ to $K$. Such a function can be written as a matrix $(v_{ij})$ with indices $i,j$ in $X$ and entries $v_{ij}$ in $K$.

The group of permutations of $X$ is denoted by $\s_X$.

Let $H$ be the Hopf $\c^*$-algebra associated to a compact quantum group of Kac type. That is, $H$ is a unital $\c^*$-algebra together with linear maps
$$\Delta :H\to H\otimes H\hskip 2cm \varepsilon :H\to\c\hskip 2cm S:H\to H$$
called comultiplication, counit and antipode and subject to the following axioms. The maps $\Delta$ and $\varepsilon$ are morphisms of $\c^*$-algebras satisfying the conditions
$$(id\otimes\Delta )\Delta =(\Delta\otimes id)\Delta\hskip 2cm (id\otimes\varepsilon )\Delta =(\varepsilon\otimes id)\Delta =id$$
and the antipode is an antimultiplicative involutive map satisfying 
$$m(id\otimes S)\Delta =m(S\otimes id)\Delta =id\hskip 2cm S^2=id$$

By \cite{wo1} the Hopf $\c^*$-algebra $H$ has a Haar functional and the dense $*$-subalgebra $H_s$ of ``smooth functions on the quantum group'' is a cosemisimple Hopf $*$-algebra.

\begin{defi}
A coaction of $H$ on $X$ is a morphism of complex $*$-algebras
$${v} :\c (X)\to \c (X)\otimes H_s$$
satisfying the coassociativity condition $({v}\otimes id){v} =(id\otimes\Delta ){v}$, the counitality condition $(id\otimes\varepsilon ){v} =id$ and the trace-preservation condition $(tr\otimes id){v} =tr (.)1$.
\end{defi}

Consider the basis of $\c (X)$ formed by the Dirac masses $\delta_i$ with $i\in X$. The coaction ${v}$ is given by the following sum.
$${v} (\delta_i)=\sum_j \delta_j\otimes v_{ji}$$

The elements $v_{ij}$ are called coefficients of the coaction. A coaction is called faithful if its coefficients generate the $*$-algebra $H$. In general we have a factorisation
$$\begin{matrix}
\c (X)&\displaystyle{\mathop{\longrightarrow}^{{v}^\prime}}&\c (X)\otimes H^\prime\cr
\ &\ &\ \cr
\ &{v}\searrow & \downarrow id\otimes\cap\cr
\ &\ &\ \cr
\ &\ & \c (X)\otimes H
\end{matrix}$$
where $H^\prime\subset H$ is the $*$-algebra generated by the coefficients of ${v}$. From the axioms for ${v}$ we get that $H^\prime$ is stable under comultiplication, so ${v}^\prime$ is a coaction as well. Thus any coaction factorises through a faithful coaction.

\medskip

As a first example, consider the Hopf $\c^*$-algebra $\c (G)$ associated to a finite group $G$. Here the comultiplication, counit and antipode are obtained by applying the $\c$ functor to the multiplication, unit and inverse map of $G$.
$$\Delta (\delta_g)=\sum_{hk=g}\delta_h\otimes\delta_k\hskip 2cm \varepsilon (f)=f(1) \hskip 2cm  S(\delta_g)=\delta_{g^{-1}}$$

If $(x,g)\mapsto gx$ is an action of $G$ on $X$, by applying the the $\c$ functor we get a map
from $\c (X)$ to $\c (X\times G)\simeq \c (X)\otimes \c (G)$, which is a coaction.
$${v} (\delta_i)=\sum_g \delta_{gi}\otimes\delta_g=\sum_j\delta_j\otimes\sum_{gi=j}\delta_g$$

This shows that ${v}$ is faithful if and only if the action map $G\to \s_X$ is injective. Thus any subgroup of the symmetric group $\s_X$ produces a faithful coaction on $X$.

\medskip

Note that the trace-preservation condition says that the action of $G$ must preserve the counting measure on $X$, and this is always true. In general this condition is not satisfied, and in case it does not hold unwanted decomposition appears in the corepresentation theory of $H$ (\cite{wa2}). In topological terms, the trace-preservation condition is needed in order to avoid ``floating circles'' in the corresponding planar algebra. See \cite{j2}, \cite{pl}.

Note also that the Kac algebra assumption $S^2=id$ is satisfied for $H=\c (G)$. In general, this assumption cannot be dropped, in the sense that for any Hopf $\c^*$-algebra coacting faithfully on $\c (X)$ the square of the antipode must be the identity.

\medskip

The axioms for ${v}$ can be translated in terms of its coefficients $v_{ij}$. First, the fact that ${v}$ is coassociative and counital translates into two ``comultiplicative'' formulae.
$$\Delta (v_{ij})=\sum_kv_{ik}\otimes v_{kj}\hskip 2cm \varepsilon (v_{ij})=\delta_{ij}$$

These conditions say that $v$ is a corepresentation of $H$ on the linear space $\c (X)$.
$$(id\otimes\Delta )v=v_{12}v_{13}\hskip 2cm (id\otimes\varepsilon )v=1$$

For the other axioms, note first that the structure of $\c (X)$ is given by the formulae
$$\delta_i^2=\delta_i\hskip 1cm\delta_i\delta_I=0\hskip 1cm 1=\sum_\iota\delta_\iota\hskip 1cm\delta_i^*=\delta_i\hskip 1cm tr(\delta_i)=1$$
valid for any $i$ and any $I\neq i$. Thus the fact that ${v}$ is a multiplicative, unital, involutive and trace-preserving translates into the following ``multiplicative'' formulae.
$$v_{ji}^2=v_{ji}\hskip 1cm v_{ji}v_{jI}=0\hskip 1cm \sum_\iota v_{j\iota}=1\hskip 1cm v_{ji}^*=v_{ji}\hskip 1cm \sum_j v_{ji}=1$$

The first four formulae say that each row of $v$ is a partition of the unity of $H$ with self-adjoint projections. In particular the matrix $v^*v$ is diagonal. The diagonal entries of $v^*v$ being the sums in the fifth formula, this formula says that $v$ satisfies $v^*v=1$. Since $v$ is a corepresentation, it follows that $v$ is unitary and that the antipode is given by the following ``third comultiplicative formula''.
$$S(v_{ij})=v_{ji}$$

By applying $S$ to the second multiplicative formula we get that the columns of $v$ are also partitions of unity with self-adjoint projections.

Summing up, in order for the map ${v}$ to be a coaction the matrix $v$ must be a corepresentation, as well as a magic biunitary in the following sense.

\begin{defi}
A matrix $v\in M_X(K)$ over a complex $*$-algebra $K$ is called a magic biunitary if its lines and columns are partitions of the unity of $K$ with projections.
\end{defi}

This vaguely reminds the definition of a magic square in combinatorics.

\medskip

Let $Z_i$ be a finite sets and let ${v}_i :\c (Z_i)\to \c (Z_i)\otimes H_i$ be coactions. Consider the corresponding magic biunitary corepresentations $v_i$. Then the matrix
$$\oplus v_i=
\begin{pmatrix}
v_1&0&\dots&0\cr
0&v_2&\dots&0\cr
\dots&\dots&\dots&\dots\cr
0&0&\dots &v_n\end{pmatrix}
$$
is a magic biunitary corepresentation of the free product $*H_i$, so we get a coaction of $*H_i$ on the disoint union $\sqcup Z_i$.
$$\oplus v_i :\c (\sqcup Z_i)\to \c (\sqcup Z_i)\otimes (*H_i)$$

This coaction is faithful if and only if all ${v}_i$ are faithful. 

As a first example, let $n_i$ be positive integers. For each $i$ the canonical action of the cyclic group $\z_{n_i}$ on itself produces a faithful coaction -- the comultiplication.
$$\Delta_i:\c (\z_{n_i})\to \c (\z_{n_i})\otimes \c (\z_{n_i})$$

The free product $*\c (\z_{n_i})$ can be computed by using the Fourier transform over cyclic groups and the distributivity of the free product of Hopf $*$-algebras with respect to the $\c^*$ functor for discrete groups (see Wang \cite{wa1}).
$$*\c (\z_{n_i})\simeq *\c^*(\z_{n_i})\simeq \c^*(*\z_{n_i})$$

We get in this way a faithful coaction
$$\oplus \Delta_i :\c (\sqcup\z_{n_i})\to \c (\sqcup\z_{n_i})\otimes \c^*(*\z_{n_i})$$

In particular for any $n\geq 4$ we can take $n_1=2$ and $n_2=n-2$ and we get a faithful coaction of $\c^*(\z_2*\z_{n-2})$ on the set with $n$ points. Thus $H$ can be infinite-dimensional. Together with the above remarks on the square of the antipode and on the trace-preservation condition, this justifies the formalism in definition 1.1.

\medskip

The universal Hopf $\c^*$-algebra coacting on $X$ is constructed by Wang in \cite{wa2}.

\begin{defi}
The universal $\c^*$-algebra $H(X)$ is defined with generators $v_{ij}$ with $i,j\in X$ and with the relations making $v$ a magic biunitary.
\end{defi}

The comultiplication, counit of antipode of $H(X)$ are constructed by using its universal property. Thus $H(X)$ is a Hopf $\c^*$-algebra coacting on $X$. For any coaction ${v} :\c (X)\to\c (X)\otimes H$ we have the following factorisation diagram.
$$\begin{matrix}
\c (X)&\displaystyle{\mathop{\longrightarrow}} &\c (X)\otimes H(X)\cr
\ &\ &\ \cr
\ &{v}\searrow &\downarrow id\otimes p\cr
\ &\ &\ \cr
\ &\ &\c (X)\otimes H
\end{matrix}$$

In particular we have a morphism $p_{com}:H(X)\to\c (G(X))$ onto the algebra of functions on the usual symmetric group $G(X)$. For sets having $1,2,3$ points this map $p_{com}$ is an isomorphism. If $X$ has at least 4 points the above free product example shows that $H(X)$ has quotients of the form $\c^*(\Gamma )$ with $\Gamma$ infinite discrete group. In particular $H(X)$ must be non commutative and infinite dimensional, and $p_{com}$ cannot be an isomorphism. See Wang \cite{wa2}.

\medskip

Some useful information about $H(X)$ can be obtained by using corepresentation theory techniques. In \cite{wo2} Woronowicz finds a general Tannaka-Krein type duality, and applies it to the algebra $H=\c (SU(N))_q$ with $q>0$. The idea is that (1) generators and relations defining $H$ translate into a presentation result of the associated tensor category, and (2) this latter equations turn to have a diagrammatic solution, and using diagrams one can compute various corepresentation theory invariants of $H$.

The same method works for $H=H(X)$. Consider the multiplication and the unit operator of the algebra $\c (X)$.
$$m:\c (X)\otimes \c(X)\to\c (X)\hskip 2cm u:\c\to\c (X)$$

These are arrows in the tensor $\c^*$-category whose objects are tensor powers of $\c (X)$, and where arrows are linear maps between such tensor powers. Let $<m,u>$ be the tensor $\c^*$-category generated by $m$ and $u$. This is by definition the smallest tensor $\c^*$-category containing $m$ and $u$. Its objects are the tensor powers of $\c (X)$, and arrows are all possible linear combinations of compositions of tensor products of maps of the form $m,u,m^*,u^*,id$, where $id$ is the identity of $\c (X)$.

Since ${v}$ is multiplicative and unital we have the following relations.
$$m\in Hom(v^{\otimes 2},v)\hskip 2cm u\in Hom(1,v)$$

Thus both $m,u$ are arrows in the category $C$ of finite dimensional corepresentations of $H(X)$. The above category $<m,u>$ can be regarded as a subcategory of $C$.

On the other hand both $m\in Hom(v^{\otimes 2},v)$ and $u\in Hom(1,v)$ give rise to a collection of relations between coefficients $v_{ij}$, which are easily seen to be equivalent to those expressing the magic biunitary condition. At the level of categories, this means that $C$ is the completion of $<m,u>$, meaning the smallest semisimple tensor $\c^*$-category containing $<m,u>$.

\medskip

This was step (1) in Woronowicz's method. What is left is step (2), namely computation of $<m,u>$ by using some combinatorial tool, such as planar diagrams. The diagrams needed here are those related to the Temperley-Lieb algebra. There are several approaches to them, beginning with those in work of Jones \cite{j1} and Kauffman \cite{k}. In this paper we use a formalism close to the one of Bisch and Jones in \cite{bj1}.

Let $\delta_X =\sqrt{\# X}$ and consider the category $TL^2(\delta_X )$ whose objects are positive integers and where arrows between $k$ and $l$ are linear combinations of Temperley-Lieb diagrams between $2k$ and $2l$ points. (Such a diagram has $k+l$ non-crossing strings joining points etc.) The operations are as follows. The composition is vertical concatenation, the tensor product is horizontal concatenation, the involution is upside-down turning and erasing a floating circle is the same as multiplying by $\delta_X$.
$$\mbox{A}\circ \mbox{B}=\begin{matrix}\mbox{B}\cr \mbox{A}\end{matrix}\hskip 2cm \mbox{A}\otimes \mbox{B}
=\mbox{A}\mbox{B}\hskip 2cm \mbox{A}^*=\forall\hskip 2cm\bigcirc =\delta_X$$

This tensor $\c^*$-category is generated by the diagrams
$$M=\delta_X^{\frac{1}{2}}\,\,\mid\cup\mid \hskip 2cm U=\delta_X^{-\frac{1}{2}}\,\,\cap$$ 
and a direct computation shows that $M\mapsto m$ and $U\mapsto u$ gives an isomorphism between $TL^2(\delta_X )$ and $<m,u>$. Thus the category of finite dimensional corepresentations of $H(X)$ is isomorphic to the completion of $TL^2(\delta_X )$. See \cite{fc}.

\medskip

Let $Y,Z$ be finite sets and consider the product $X=Y\times Z$. The projection $X\to Y$ gives a unital embedding of $\c (Y)$ into $\c (X)$, so we can consider coactions on $X$ which restrict to coactions on $Y$.
$$\begin{matrix}
\c (X)&\displaystyle{\mathop{\longrightarrow}^{{v}}} &\c (X)\otimes H\cr
\ &\ &\ \cr
\cup &\ &\cup\otimes id\cr
\ &\ &\ \cr
\c (Y) &\displaystyle{\mathop{\longrightarrow}} &\c (Y)\otimes H
\end{matrix}$$

For a coaction ${v}$ to restrict to $\c (Y)$, it must send $\c (Y)$ into $\c (Y)\otimes H$.

\begin{defi}
The Hopf $\c^*$-algebra $H(X\to Y)$ is the quotient of $H(X)$ by the closed two-sided $*$-ideal generated by the relations coming from the equality $ev=ve$, where $e:\c (X)\to\c (X)$ is the linear map implemented by the projection $X\to Y$.
\end{defi}

The Hopf $\c^*$-algebra $H(X\to Y)$ coacts on $X$ and the coaction restricts to $Y$. For any coaction on $X$ which restricts to $Y$ we have the following factorisation diagram.
$$\begin{matrix}
\c (X)&\displaystyle{\mathop{\longrightarrow}} &\c (X)\otimes H(X\to Y)\cr
\ &\ &\ \cr
\ &{v}\searrow &\downarrow id\otimes p\cr
\ &\ &\ \cr
\ &\ &\c (X)\otimes H
\end{matrix}$$

The category of corepresentations of $H(X\to Y)$ can be computed by using the same techniques as for $H(X)$, with the Fuss-Catalan diagrams of Bisch and Jones \cite{bj1}.

First, the universal property of $H(X\to Y)$ translates into the fact that its category of corepresentations is the completion of the tensor $\c^*$-category generated by $m$, $u$ and $e$. Let $\delta_Y=\sqrt{\# Y}$ and $\delta_Z=\sqrt{\# Z}$ and consider the category $FC(\delta_Y,\delta_Z)$ whose objects are positive integers and where the arrows between $l$ and $k$ are linear combinations of Fuss-Catalan diagrams between $4l$ points and $4k$ points. These are by definition Temperley-Lieb diagrams with both rows of points colored
$$y,z,z,y,y,z,z,y,\ldots$$
such that strings join points having the same color. Then $FC(\delta_Y,\delta_Z)$ is a tensor $\c^*$-category, with operations given by usual horizontal and vertical concatenation and upside-down turning of diagrams, plus the following rule: erasing a $y$-colored or $z$-colored circle is the same as multiplying the diagram by $\delta_Y$ or $\delta_Z$.
$$\mbox{\tiny{y-colored}}\rightarrow\bigcirc =\delta_Y\hskip 2cm
\mbox{\tiny{z-colored}}\rightarrow\bigcirc =\delta_Z$$

With $\delta_X=\delta_Y\delta_Z=\sqrt{\# X}$ this category is generated by the diagrams
$$M=\delta_X^{\frac{1}{2}}\,\,
\mid\mid\!\Cup\!\mid\mid\hskip 2cm U=\delta_X^{-\frac{1}{2}}\,\, \Cap\hskip 2cm E
=\delta_Z^{-\frac{1}{2}}\,\,\mid\begin{matrix}\cup\cr\cap\end{matrix}\mid$$
and computation shows that $M\mapsto m$, $U\mapsto u$ and $E\mapsto e$ gives an isomorphism between $FC(\delta_Y,\delta_Z)$ and $<m,u,e>$. Thus the category of finite dimensional corepresentations of $H(X\to Y)$ is isomorphic to the completion of $FC(\delta_Y,\delta_Z)$. See \cite{fc}.

\medskip

Finally, let us mention that quite exotic examples, not to be discussed in this paper, can be obtained from magic biunitary matrices. Let $X,Y$ be finite sets and let $(u_{ij})$ be a magic biunitary with $i,j\in X$ and $u_{ij}\in M_Y(\c )$. The universal property of $H(X)$ gives a representation on $M_Y(\c)$ sending the fundamental corepresentation $v$ to $u$.
$$\pi :H(X)\to M_Y(\c)\hskip2cm (id\otimes\pi )v=u$$

By dividing $H(X)$ by a suitable ideal we get a factorisation $\pi =\pi_up_u$ through a minimal possible Hopf $\c^*$-algebra $H(u)$. This coacts on $\c (X)$ via $v_u=(id\otimes p_u)v$.
$$\pi_u :H(u)\to M_Y(\c)\hskip2cm (id\otimes\pi_u )v_u=u$$

For instance if $G$ is a subgroup of $\s_Y$ and $u$ is the magic biunitary corresponding to the action of $G$ on $Y$ then $H(u)=\c (G)$ and $v_u=u$.

More complicated examples come from lattice models and commuting squares. For instance complex Hadamard matrices and Krishnan-Sunder biunitary permutation matrices are known to produce magic biunitaries (cf. \cite{js}, \cite{bh2}, \cite{sq}). So far nothing seems to be known about the associated Hopf $\c^*$-algebras, besides their existence.

\section{Coactions on finite metric spaces}

Let $(X,d)$ be a finite metric space. That is, we are given a finite set $X$ and a function $d:X\times X\to\r_+$ which is symmetric, which has 0 on the diagonal and positive values outside the diagonal, and whose entries satisfy the triangle inequality.

Let $G$ be a group acting on $X$. Consider the corresponding coaction of $\c (G)$ on $X$.
$${v} :\c (X)\to\c (X)\otimes \c (G)\hskip 2cm {v} (\delta_i)=\sum_g \delta_{gi}\otimes\delta_g$$

The products $dv$ and $vd$ are given by the following formulae.
$$(dv)_{ij}=\sum_g d(i,gj)\delta_g\hskip 2cm (vd)_{ij}=\sum_gd(g^{-1}i,j)\delta_g$$

Thus the coaction of $\c (G)$ on $X$ preserves the metric if and only if the action of $G$ on $X$ preserves the metric. This suggests the following definition.

\begin{defi}
Let $H$ be the Hopf $\c^*$-algebra associated to a compact quantum group of Kac type. A coaction of $H$ on $X$
$${v} :\c (X)\to\c (X)\otimes H\hskip 2cm {v}(\delta_i)=\sum_j\delta_j\otimes v_{ji}$$
is said to preserve the metric if the matrix  $v=(v_{ij})$ commutes with $d=(d(i,j))$.
\end{defi}

Let $H(X,d)$ be the quotient of $H(X)$ by the relations coming from $dv=vd$. Then $H(X,d)$ coacts on $\c (X)$ and preserves the metric. For any other coaction which preserves the metric we have the following factorisation diagram.
$$\begin{matrix}
\c (X)&\displaystyle{\mathop{\longrightarrow}} &\c (X)\otimes H(X,d)\cr
\ &\ &\ \cr
\ &{v}\searrow &\downarrow id\otimes p\cr
\ &\ &\ \cr
\ &\ &\c (X)\otimes H
\end{matrix}$$

In particular we have the following factorisation of the universal coaction
$$\begin{matrix}
\c (X)&\displaystyle{\mathop{\longrightarrow}} &\c (X)\otimes H(X,d)\cr
\ &\ &\ \cr
\ &{v}_{com}\searrow &\downarrow id\otimes p_{com}\cr
\ &\ &\ \cr
\ &\ &\c (X)\otimes \c (G(X,d))
\end{matrix}$$
where $G(X,d)$ is the automorphism group of $(X,d)$ and ${v}_{com}$ is the corresponding coaction. We will see that in many cases the projection $p_{com}$ is an isomorphism.

The simplest example of a finite metric space is a finite set with the 0--1 distance on it. More generally, let us call simplex a finite metric space having the same distance between all pairs of disjoint points. 

\begin{prop}
If $(X,d)$ is a simplex then any coaction on $X$ preserves the metric. We have $H(X,d)=H(X)$.
\end{prop}

\begin{proof}
The distance on a simplex is given by the formula
$$d=a (\one -1)$$
where $a$ is the non-zero value of the distance, $1$ is the identity matrix and $\one$ is the matrix having 1 everywhere.
$$\one =
\begin{pmatrix}
1&\dots&1\cr
\dots&\dots&\dots\cr
1&\dots&1
\end{pmatrix}$$

Let $v$ be the magic biunitary corepresentation corresponding to a coaction on $X$. Since right and left multiplication by $\one$ is making sums on rows and columns, $\one$ commutes with $v$. Thus $v$ commutes with $d$ as well.
\end{proof}

The metric analogue of the free product example is as follows. Let $(Y,D)$ be a metric space and for each $i\in Y$ let $(Z_i,d_i)$ be a metric space. Assume that the following condition is satisfied.
$$2D(i,j)\geq d_i(p,q)$$

This must hold for any $i,j\in Y$ and for any $p,q\in Z_i$. The disjoint union $X=\sqcup Z_i$ becomes a metric space, with distance $d$ equal to $d_i$ on each $Z_i$ and equal to $D(i,j)$ between points of $Z_i$ and $Z_j$ for $i\neq j$. It is useful to keep in mind the following picture: $(X,d)$ is made of a metric ``macrospace'' $(Y,D)$ having around each of its points $i$ a metric ``microspace''  $(Z_i,d_i)$.

\begin{prop}
If ${v}_i :\c (Z_i)\to\c (Z_i)\otimes H_i$ are coactions then the coaction
$$\oplus v_i :\c (X)\to \c (X)\otimes (*H_i)$$
preserves the metric $d$ if and only if each $v_i$ preserves the metric $d_i$.
\end{prop}

\begin{proof}
If points of $Y$ are labelled from 1 to $n$ then the distance on $X$ is
$$d=\begin{pmatrix}
d_1&D(1,2)\one &\dots& D(1,n)\one\cr 
D(2,1)\one &d_2&\dots& D(2,n)\one\cr 
\dots&\dots&\dots&\dots\cr 
D(n,1)\one &D(n,2)\one&\dots& d_n
\end{pmatrix}$$
where $\one$ denotes rectangular matrices of suitable size filled with $1$'s. Since multiplication by $\one$ is making sums on rows or columns, $\one$ commutes with each $v_i$. Thus $\oplus v_i$ commutes with $d$ is and only if each $v_i$ commutes with $d_i$.
\end{proof}

As an example, let $Y$ be the segment of length $1$ and let $Z_i$ be two segments of length $\sqrt{2}$. The metric space $X$ is a square. Since $\z_2$ acts on both segments, now diagonals of the square, it follows that $\c^*(\z_2*\z_2 )$ coacts on the square.

More generally, let $Y$ be a simplex with nonzero distance $A$ and $Z$ be a simplex with nonzero distance $a$ and assume that $A\neq a$ and $2A\geq a$.  For any $i\in Y$ define $Z_i=Z$. The disjoint union in the above example is $X=Y\times Z$, with the following metric.
$$d=\begin{pmatrix}
a(\one -1)&A\one &\dots& A\one\cr 
A\one &a(\one -1)&\dots& A\one\cr 
\dots&\dots&\dots&\dots\cr 
A\one &A\one&\dots& a(\one -1)
\end{pmatrix}$$

It is useful to keep in mind the following picture: $X$ is a space made of a ``macrosimplex'' $Y$, having around each of its points a copy of the ``microsimplex''  $Z$.

\begin{prop}
A coaction on $X$ preserves the metric if and only if it restricts to $Y$. We have $H(X,d)=H(X\to Y)$.
\end{prop}

\begin{proof}
We have the following formula
$$d=A\one + (a-A)(\# Z)e-a1$$
where $1$ is the identity matrix and where $e$ is the following operator.
$$e=\frac{1}{\# Z}\begin{pmatrix}
\one&0&\dots&0\cr
0& \one&\dots&0\cr
\dots&\dots&\dots&\dots\cr
0&0&\dots &\one
\end{pmatrix}$$

Let $v$ be a coaction on $X$. Since commutation of $v$ with $1$ and with $\one$ is automatic, $v$ preserves the distance if and only if it commutes with $e$. But $e$ is the projection onto the subalgebra $\c(Y)$ of $\c (X)$.
\end{proof}

Given a finite metric space $(X,d)$ and a group $G$ acting on $X$, the fact that $G$ preserves or not the metric only depends on some equivalence class of the metric: the one which takes into account the equalities between various pairs of distances, but neglects their precise values. In other words, the induced action of $G$ on $X\times X$ must preserve the 0--1 matrix of Kronecker symbols
$$K(i,j,k,l)=\delta_{d(i,j),d(k,l)}$$

Instead of using $K$ or other purely combinatorial objects, it is convenient to use the metric and to draw pictures, but to keep in mind that the set of values of the distance function is irrelevant. We will call this values ``colors'' and the picture of a metric space will consist of points plus colored lines between them, with the convention that same color means same distance. Because of typographical reasons colors in this paper will be different kinds of lines, plus a ``void color'' that is not pictured. The void color is used in order to simplify pictures and will be typically a most frequent color, i.e. a number $a$ which maximizes the following quantity.
$$\# \{ (i,j)\vert d(i,j)=a\}$$

As an example, here is a rectangle which is not a square. The void color is here the length of the diagonals, which is different from the lengths of the sides.
$$\vert\begin{matrix}\cdots\cdots\cr\cdots\cdots\end{matrix}\vert$$

Let $(X,d)$ be a metric space with color set $C\subset\r_+$. For any color $a\in C$ we denote by $d_a$ the 0--1 matrix having $1$ if $d(i,j)=a$ and 0 elsewhere. The following formula will be called the color decomposition of the metric.
$$d=\sum_{a\in C}a d_a$$

For instance the color decomposition for the rectangle having sides of lengths $a$ and $c$ and diagonals of lengths $b$ is as follows (in order to obtain usual $4\times 4$ matrices, the numers 1,2,3,4 are assigned clockwise to the vertices of the rectangle).
$$\begin{pmatrix}0&a&b&c\cr a&0&c&b\cr b&c&0&a\cr c&b&a&0\end{pmatrix}=
a\begin{pmatrix}0&1&0&0\cr 1&0&0&0\cr 0&0&0&1\cr 0&0&1&0\end{pmatrix}+
b\begin{pmatrix}0&0&1&0\cr 0&0&0&1\cr 1&0&0&0\cr 0&1&0&0\end{pmatrix}+
c\begin{pmatrix}0&0&0&1\cr 0&0&1&0\cr 0&1&0&0\cr 1&0&0&0\end{pmatrix}$$

The equivalence classes of the spaces $(X,d_a )$ will be called color components of $(X,d)$. Note that each $(X,d_a)$ is not exactly a metric space, because the ``distances'' $d_a$ can have 0 values outside the diagonal. But each such space corresponds to a unique equivalence class of a metric space, which can be obtained by replacing all non-diagonal $0$ entries of $d_a$ with some positive number which is different from 1, bigger than 0.5 and smaller than 2. For instance each color component of the rectangle is a square.

With these notations, the equivalence class of a metric space $(X,d)$ is given by its color components. The above remark about actions of groups on equivalence classes has the following reformulation. A group $G$ acting on $X$ preserves the metric if and only if preserves the metric of all color components of $(X,d)$.

This has the following analogue.

\begin{prop}
A coaction of $H$ on $X$ preserves the metric if and only if it preserves all its color components.
\end{prop}

\begin{proof}
Assume that $H$ preserves the metric and consider the multiplication and comultiplication of $\c (X)$.
$$m:\delta_i\otimes \delta_j\mapsto \delta_i\delta_j\hskip 2cm c :\delta_i\mapsto \delta_i\otimes \delta_i$$

Then $m$ and $c$ intertwine $v^{\otimes 2}$ and $v$. It follows that their iterations
$$m^{(k)}:\delta_1\otimes\ldots\otimes \delta_k\mapsto \delta_1\ldots \delta_k\hskip 2cm c^{(k)} :\delta_i\mapsto \delta_i\otimes\ldots\otimes \delta_i$$
are in $Hom(v^{\otimes k},v)$ and $Hom(v,v^{\otimes k})$ respectively. In particular the operator
$$d^{(k)}=m^{(k)}d^{\otimes k}c^{(k)}:\delta_i\mapsto\sum_j d(i,j)^k\delta_j$$
is in the commutant $End(v)$ of $v$. Since $d^{(k)}$ is obtained from $d$ by raising all its entries to the power $k$, it has the following decomposition.
$$d^{(k)}=\sum_{a\in C}a^k d_a$$

If $A$ is the biggest color then the sequence $A^{-k}d^{(k)}$ of elements of $End(v)$ converges to $d_A$, so $d_A$ has to be an element of $End(v)$. By substracting $Ad_A$ from $d$ we get in this way by induction that each $d_a$ is in $End(v)$.
\end{proof}

\begin{coro}
Let $v:\c (X)\to \c (X)\otimes H$ be a metric-preserving coaction. If $K$ is an eigenspace of some color component $d_a$ then $v(K)\subset K\otimes H$. 
\end{coro}

\begin{proof}
Consider the orthogonal projection $P_K$ onto $K$. Since $K$ is an eigenspace of $d_a$ this projection is in the $\c^*$-algebra generated by $d_a$ in $\l (\c (X))$. Thus from $d_a\in End(v)$ we get that $P_K$ is in $End(v)$ as well. On the other hand we have $v(K)\subset K\otimes H$ if and only if $v$ commutes with $P_K$.
\end{proof}

\section{Rules for computing the Hopf algebra}

In this section we find some rules for the computation of $H(X,d)$ and its corepresentations. There is a lot of structure on finite metric spaces that can be exploited. For instance simple geometric properties of $(X,d)$ lead to invariant subspaces or subalgebras of $v$ of dimension 1,2 or 3 and the universal properties of the Hopf $\c^*$-algebras $\c^*(\z)$, $\c (SU(2))$, $\c (SU(2))_{-1}$, $\c (\z_2)$ and $\c (\z_3)$ can be used. However, the computation of $H(X,d)$ seems to be hard in the general case and the set of needed ``rules'' might be infinite. In fact we don't know what happens in the $n=8$ case and our less ambitious project here is just to find enough rules for dealing with the $n\leq 7$ case.

We begin with reformulations of three facts from the previous sections.

\begin{lemm}[Triangle Rule]
If $X$ has at most 3 points then $H(X,d)=\c (G(X,d))$.
\end{lemm}

\begin{lemm}[Simplex Rule]
If $X$ is a simplex having at least 4 points the category of finite dimensional corepresentations of $H(X,d)$ is the completion of the category of Temperley-Lieb diagrams of parameter $\delta_X = \sqrt{\#X}$.
\end{lemm}

\begin{lemm}[Duplex Rule]
If $X=Y\times Z$ is a product of simplices then the category of finite dimensional corepresentations of $H(X,d)$ is the completion of the category of Fuss-Catalan diagrams of parameters $\delta_Y=\sqrt{\#Y}$ and $\delta_Z=\sqrt{\# Z}$.
\end{lemm}

Let $(X,d)$ be a finite metric space with color set $C\subset\r_+$. For any partition
$$C=\{ a_1,\ldots ,a_m\}\cup \{ b_1,\ldots ,b_n\}\cup\{ c_1,\ldots ,c_p\}\cup\ldots$$
we get a new space $(X^\prime ,d^\prime )$ by identifying all $a$-colors, all $b$-colors, all $c$-colors etc. This is not exactly a metric space, but it corresponds to a unique equivalence class of a metric space (cf. comments in 2.). We say that $(X^\prime ,d^\prime )$ is a decoloration of $(X,d)$.

As an example, consider the rectangle. Its color set has three elements $\{ a,b,c\}$ so there are 5 possible partitions of it. The corresponding decolorations are the rectangle itself (no decoloration), the square (appearing 3 times, obtained by identifying 2 colors) and the 4-simplex (obtained by identifying all colors).

Each color component of $(X,d)$ being a decoloration of $(X,d)$, proposition 2.4 has the following slight generalisation plus reformulation.

\begin{lemm}[Decoloration Rule]
If $(X^\prime ,d^\prime )$ is a decoloration of $(X,d)$ then $H(X,d)$ is a quotient of $H(X^\prime ,d^\prime )$.
\end{lemm}

This can be used as follows. Let $(X,d)$ be a finite metric space and assume that we have found a decoloration $(X^\prime ,d^\prime )$ such that $H(X^\prime ,d^\prime )$ is commutative. Then $H(X,d)$ must be commutative as well, so it has is equal to $\c (G(X,d))$. A similar corollary holds with ``commutativity'' replaced by the ``transitivity'' notion discussed below.

We say that a color $a$ is cyclic if there exists a labelling of the points of $X$ such that the $a$ color connects 1 to 2 to 3 etc.

\begin{lemm}[Cycle Rule]
Assume there exists a cyclic color and that the space has $n\geq 5$ points. Then $H(X,d)=\c (G(X,d))$.
\end{lemm}

\begin{proof}
By using the decoloration rule it is enough to prove that if a Hopf $\c^*$-algebra $H$ coacts faithfully on the cycle $(X,d_a)$ then $H$ is commutative.

The matrix of the cycle is
$$M=\begin{pmatrix}
0     &1     &0     &0     &\ldots&0     &0     &1     \cr
1     &0     &1     &0     &\ldots&0     &0     &0     \cr
0     &1     &0     &1     &\ldots&0     &0     &0     \cr
\ldots&\ldots&\ldots&\ldots&\ldots&\ldots&\ldots&\ldots\cr
0     &0     &0     &0     &\ldots&1     &0     &1     \cr
1     &0     &0     &0     &\ldots&0     &1     &0
\end{pmatrix}$$

Let $w$ be a $n$-th root of 1 and consider the vector
$$f=\begin{pmatrix}1\cr w\cr w^2\cr\dots\cr w^{n-1}\end{pmatrix}$$

Then $f$ is an eigenvector of $M$ with eigenvalue $w+w^{n-1}$. Now by taking $w$ to be a primitive $n$-th root of 1, it follows that all vectors $1,f,f^2,\ldots ,f^{n-1}$ are eigenvectors of $M$. The invariant subspaces of $M$ are
$$\c 1,\, \c f\oplus\c f^{n-1},\, \c f^2\oplus\c f^{n-2},\ldots$$
where the last subspace has dimension 1 or 2 depending on the parity of $n$. Write
$${v} (f)=f\otimes a+f^{n-1}\otimes b$$
with $a,b\in H$ (cf. corollary 2.1, to be used many times in what follows). By taking the square of this equality we get
$${v} (f^2)=f^2\otimes a^2+f^{n-2}\otimes b^2+1\otimes (ab+ba)$$

It follows that $ab=-ba$ and that ${v} (f^2)$ is given by the formula
$${v} (f^2)=f^2\otimes a^2+f^{n-2}\otimes b^2$$

By multiplying this with ${v} (f)$ we get
$${v} (f^3)=f^3\otimes a^3+f^{n-3}\otimes b^3+f^{n-1}\otimes ab^2+f\otimes ba^2$$

But $n\geq 5$ implies that 1 and $n-1$ are different from 3 and from $n-3$, so we must have $ab^2=ba^2=0$. Now by using this and $ab=-ba$ we get by induction on $k$ that
$${v} (f^k)=f^k\otimes a^k+f^{n-k}\otimes b^k$$
for any $k$. In particular for $k=n-1$ we get
$${v} (f^{n-1})=f^{n-1}\otimes a^{n-1}+f\otimes b^{n-1}$$

On the other hand we have $f^*=f^{n-1}$, so by applying $*$ to ${v} (f)$ we get
$${v} (f^{n-1})=f^{n-1}\otimes a^*+f\otimes b^*$$

Thus $a^*=a^{n-1}$ and $b^*=b^{n-1}$. Together with $ab^2=0$ this gives
$$(ab)(ab)^*=abb^*a^*=ab^na^{n-1}=(ab^2)b^{n-2}a^{n-1}=0$$
and from positivity we get $ab=0$. Together with $ab=-ba$ this shows that $a$ and $b$ commute. On the other hand $H$ is generated by the coefficients of ${v}$, which are powers of $a$ and $b$, so $H$ must be commutative.
\end{proof}

Assume that the metric space $(X,d)$ has an even number of points. A star in $X$ is given by two colors $a$ and $b$ and a labelling of points of $X$ from $1$ to $2k$ such that $(X,d_a)$ consists of the segments $(1,2)$, $(3,4)$,\ldots ,$(2k-1,2k)$ and $(X,d_b )$ consists of the segments $(2,3)$, $(4,5)$,\ldots ,$(2k-2,2k-1)$. (If points of $X$ are arranged clockwise this doesn't quite look like a star; but a suitable rearrangement of points justifies terminology.)

\begin{lemm}[Star Rule]
If $X$ has a star then $H(X,d)=\c (G(X,d))$.
\end{lemm}

\begin{proof}
By using the projection $p_{com}$ it is enough to prove that $H(X,d)$ is commutative.

Let $(X^\prime ,d^\prime )$ be the decoloration of $(X,d)$ obtained by identifying the colors $a$ and $b$. The star in $X$ becomes a cycle in $X^\prime$. From the decoloration rule we get that $H(X,d)$ is a quotient of $H(X^\prime,d^\prime)$, so if $H(X^\prime,d^\prime)$ is commutative then $H(X,d)$ is commutative as well. This is the case for $k\geq 3$, because we have $2k\geq 6$ and the cycle rule applies.

In the remaining $k=2$ case $X$ must be a rectangle. Consider the vectors
$$1=\begin{pmatrix}1\cr 1\cr 1\cr 1\end{pmatrix}\hskip 1cm
A=\begin{pmatrix}1\cr -1\cr 1\cr -1\end{pmatrix}\hskip 1cm
B=\begin{pmatrix}1\cr 1\cr -1\cr -1\end{pmatrix}\hskip 1cm
C=\begin{pmatrix}1\cr -1\cr -1\cr 1\end{pmatrix}$$

By using the color decomposition of the rectangle given in 3. we see that these are eigenvectors for all three matrix distances, corresponding to different eigenvalues. By corollary 2.1 the three 1-dimensional spaces spanned by $A$, $B$ and $C$ must be left invariant by the coaction, so we must have
$${v} (A)=A\otimes a\hskip 1cm{v} (B)=B\otimes b\hskip 1cm{v} (C)=C\otimes c$$
for some elements $a,b,c\in H$. We have the following relations between $A,B,C$
$$AB=BA=C\hskip 1cm AC=CA=B\hskip 1cm CB=BC=A$$
and by applying ${v}$ we get from the multiplicativity of ${v}$ that $a,b,c$ satisfy similar formulae, hence commute. But $a,b,c$ generate $H$, so $H$ is commutative.
\end{proof}

A coaction ${v} :\c (X)\to \c (X)\otimes H$ is called transitive if ${v} (f)=f\otimes 1$ implies that $f$ is a constant function. This is an extension of the classical notion of transitivity: in the $H=\c (G)$ case we have the formulae
$${v} (f)=\sum_j\delta_j\otimes\sum_gf(g^{-1}j)\delta_g\hskip 2cm f\otimes 1=
\sum_j\delta_j\otimes\sum_gf(j)\delta_g$$
so the coaction of $\c (G)$ is transitive if and only if the action of $G$ has one orbit.

A metric space $(X,d)$ is called homogenous if its automorphism group acts transitively, meaning has one orbit. This is the same as asking that for any $x,y\in X$ there is a metric-preserving permutation $g\in\s_X$ such that $gx=y$.

One can say that $(X,d)$ is ``quantum homogenous'' is the coaction of $H(X,d)$ is transitive in the above sense. We have here the following necessary condition.

\begin{lemm}[Magic Rule]
If the coaction of $H(X,d)$ is transitive then the following condition must be satisfied. For any two points $x,y\in X$ there exists a bijection $\sigma :X\to X$ such that $\sigma (x)=y$ and such that $d(x,z)=d(y,\sigma (z))$ for any $z$.
\end{lemm}

The terminology comes from the following fact. Assume that $X$ has $n$ points and that there are exactly $n-1$ colors. By replacing these colors with the numbers from 1 to $n-1$, the matrix distance must be a magic matrix, i.e. a $n\times n$ matrix having the numbers from 0 to $n-1$ on each row and column. As an example, here is the magic matrix corresponding to the rectangle.
$$\begin{pmatrix}0&1&2&3\cr 1&0&3&2\cr 2&3&0&1\cr 3&2&1&0\end{pmatrix}$$

\begin{proof}
For any color $a$ we have $d_a\in End(v)$, so $d_a (1)$ is fixed by the coaction. But $d_a(1)$ is scalar if and only if the matrix $d_a$ has the same number of 1's on each row.
\end{proof}

We say that a color $a$ is bicyclic if $(X,d_a)$ consists of two disjoint cycles $C_n\sqcup C_N$.

\begin{lemm}[Bicycle Rule]
Assume there exists a bicyclic color component of the form $C_n\sqcup C_N$ with $n\neq N$. Then the coaction of $H(X,d)$ is not transitive.
\end{lemm}

\begin{proof}
By using the decoloration rule it is enough to prove that if a Hopf $\c^*$-algebra $H$ coacts on the bicycle then the coaction is not transitive.

If one of $n$ and $N$ is $2$ then the magic rule shows that $H(X,d)$ is not transitive. So we may assume $n,N\geq 3$. Label the points of $X$ such that the matrix of distances is
$$M=\begin{pmatrix}M_n&0\cr 0&M_N\end{pmatrix}$$
where $M_n$ and $M_N$ are the matrices of the two cycles. Let $w$ be a primitive $n$-th root of unity and consider the vector $f$ consisting of powers of $w$. Then the eigenvalues and eigenvectors of the matrix $(^{M_n}_0{\ }^0_0)$ are given by
$$2:\left[\e{0}{1}\right] ,\, (w+w^{n-1}):\left[\e{0}{f},\e{0}{f^{n-1}}\right] , (w^2+w^{n-2}):\left[\e{0}{f^2},\e{0}{f^{n-2}}\right]\,\ldots$$
where each eigenvalue is followed by its eigenvectors and where the end of the sequence depends on the parity of $n$. Similarly, the matrix $(^0_0{\ }^0_{M_N})$ gives
$$2:\left[\e{1}{0}\right] ,\, (W+W^{N-1}):\left[\e{F}{0},\e{F^{N-1}}{0}\right] , (W^2+W^{N-2}):\left[\e{F^2}{0},\e{F^{N-2}}{0}\right]\,\ldots$$
where $W$ is a primitive $N$-th root of unity and $F$ is the vector consisting of powers of $W$. Since $n\neq N$ it follows that the decomposition corresponding to $M$ is
$$\left[\e{0}{1},\e{1}{0}\right], \,
\begin{matrix} 
\left[\e{0}{f},\e{0}{f^{n-1}}\right],
\left[\e{0}{f^2},\e{0}{f^{n-2}}\right],\,\ldots\cr
&\cr
\left[\e{F}{0},\e{F^{N-1}}{0}\right],
\left[\e{F^2}{0},\e{F^{N-2}}{0}\right],
\,\ldots
\end{matrix}$$
where the two series might interfere on their last entries, depending on the parity of $n$ and $N$. Consider the first subspace, corresponding to the eigenvalue 2.
$$\left[\e{0}{1},\e{1}{0}\right] =\left[\e{1}{1},\e{-1}{1}\right]$$

By corollary 2.1 this space is invariant by ${v}$, so we can find $\xi, u\in H$ such that
$${v}\e{-1}{1}=\e{1}{1}\otimes\xi +\e{-1}{1}\otimes u$$

From the trace-preservation property of ${v}$ we get that $tr\e{-1}{1}1=0$ is equal to
$$(tr\otimes id){v} \e{-1}{1}=tr\e{1}{1}\xi +tr\e{-1}{1}u=2\xi$$

It follows that $\xi =0$, so we have
$${v}\e{-1}{1}=\e{-1}{1}\otimes u$$

From this we get that ${v}$ is given on the subspace by
$${v}\e{0}{1}=\e{0}{1}\otimes\frac{1+u}{2}+\e{1}{0}\otimes\frac{1-u}{2}$$
$${v}\e{1}{0}=\e{0}{1}\otimes\frac{1-u}{2}+\e{1}{0}\otimes\frac{1+u}{2}$$

Assume that one of $n$ and $N$ is greater than 5, say $n\geq 5$. By arguing like in the proof of the cycle rule we get from a formula of type
$${v} \e{0}{f}=\e{0}{f}\otimes a+\e{0}{f^{n-1}}\otimes b$$
that the following formula is true for any $k$.
$${v} \e{0}{f^k}=\e{0}{f^k}\otimes a^k+\e{0}{f^{n-k}}\otimes b^k$$

In particular for $k=n$ we get
$${v} \e{0}{1}=\e{0}{1}\otimes a^n+\e{0}{1}\otimes b^n=\e{0}{1}\otimes (a^n+b^n)$$
and comparison with the above formula for ${v} \e{0}{1}$ gives $u=1$. It follows that $\e{-1}{1}$ is a fixed vector of the coaction, so the coaction is not transitive.

The remaining case is $n,N\leq 4$. Since $n,N$ are different and both are greater than 3, one is 3 and the other is 4. Assume $N=4$ and consider the vector
$$\gamma =F^2=\begin{pmatrix}1\cr-1\cr1\cr-1\end{pmatrix}$$

This is the unique eigenvector of $M_N$ corresponding to the eigenvalue $0$. Since $n=3$ is odd, $M_n$ has no eigenvalue 0. It follows that $\e{\gamma}{0}$ is the unique eigenvector of $M$ corresponding to the eigenvalue $0$, so we must have
$${v}\e{\gamma}{0}=\e{\gamma}{0}\otimes c$$
for some $c\in H$. By taking the square of this equality we get
$${v}\e{1}{0}=\e{1}{0}\otimes c^2$$
and comparison with the above formula for ${v} \e{0}{1}$ gives $u=1$. Once again, it follows that $\e{-1}{1}$ is a fixed vector of the coaction, so the coaction is not transitive.
\end{proof}

\section{Classification for $n\leq 7$}

The purpose of this section is to compute the corepresentation theory invariants of the Hopf algebra $H(X,d)$, in case the metric space $(X,d)$ has at most 7 points.

In order to cut off unwanted complexity, there is a natural condition to be put on $(X,d)$, namely quantum homogeneity. Recall from previous section that $(X,d)$ is said to be quantum homogenous when the universal coaction $v:\c (X)\to \c(X)\otimes H(X,d)$ is transitive, meaning that $v(f)=f\otimes 1$ implies that $f$ is a scalar function.

The motivating remark here is that functions satisfying $v(f)=f\otimes 1$ are those belonging to the relative commutant of the associated subfactor (see \cite{pl}), so quantum homogeneity corresponds to the well-known irreducibility condition in subfactors.

It is clear from definitions that a homogenous space is quantum homogenous. It turns out that for spaces having at most 7 points the converse holds. We do not know if these notions are equivalent in general.

\medskip

The decoloration rule shows that all color components of a quantum homogenous space must be quantum homogenous. These color components are graphs, so our first task is to prove that for a graph having at most 7 vertices, homogeneity is equivalent to quantum homogeneity.

The complement of a graph $(X,d)$ with vertex set $X$ and incidency matrix $d$ is by definition the graph $(X,\one -1-d)$, where $\one$ is the matrix filled with 1's. Since $v\one=\one v=\one$ for any magic biunitary matrix $v$, we have $H(X,d)=H(X,\one -1-d)$.

\begin{lemm}
A graph with $n\leq 7$ vertices is homogenous if and only if it is quantum homogenous. There are $24$ such graphs.

(1) point.

(2) $2$ points, segment.

(3) $3$ points, triangle.

(4) $4$ points, tetrahedron, $2$ segments, square.

(5) $5$ points, $5$-simplex, $5$-cycle.

(6) $6$ points, $6$-simplex, plus $6$-cycle, $3$ segments, $2$ triangles, and their complements.

(7) $7$ points, $7$-simplex, $7$-cycle and its complement.
\end{lemm}

\begin{proof}
Homogenous implies quantum homogenous, and all graphs in the list are homogenous. It remains to prove that any quantum homogenous graph with $n\leq 7$ vertices is in the list.

Let $(X,d)$ be such a graph. The magic rule applies and shows that the incidency matrix $d$ has the same number of 1's, say $k$, on each of its rows and columns. Thus $(X,d)$ must be a $k$-regular graph, for a certain number $k$. By replacing the graph with its complement we may assume that that the valence $k$ of the graph is smaller or equal than the valence $n-k-1$ of its complement.

Summing up, it remains to prove the following assertion. Any quantum homogenous $k$-regular graph with $n$ vertices, $2k+1\leq n\leq 7$, is in the list.

For $k=0$ the graph must be $n$ points and no edges, in the list.

For $k=1$ the only graphs are the 2 or 3 segments, in the list.

For $k=2$ the graph must be a union of $m$-cycles with values of $m$ greater than 3. The bicycle rule shows that the $m$'s must be the same for all $m$-cycles, and the graphs left are the $n$-cycles with $3\leq n\leq 7$, plus the 2 triangles, all of them in the list.

For $k\geq 3$ we must have $n\geq 2\times 3+1=7$. The only possible case is $k=3$ and $n=7$, but this is excluded by the observation that the incidency matrix must have $nk/2$ values of $1$ above the diagonal.
\end{proof}

\begin{theo}
Let $(X,d)$ be a metric space having at most 7 points. Then one of the following happens.

(1) The coaction of $H(X,d)$ is not transitive.

(2) $H(X,d)$ is the algebra of functions on the usual automorphism group of $(X,d)$.

(3) $X$ is a simplex and the category of finite dimensional corepresentations of $H(X,d)$ is the completion of the category of Temperley-Lieb diagrams of parameter $\delta_X=\sqrt{\#X}$.

(4) $X=Y\times Z$ is a product of simplices and the category of finite dimensional corepresentations of $H(X,d)$ is the completion of the category of Fuss-Catalan diagrams of parameters $\delta_Y=\sqrt{\#Y}$ and $\delta_Z=\sqrt{\# Z}$.
\end{theo}

\begin{proof}
If $X$ is a simplex the triangle or simplex rule applies and shows that we are in situation (2) or (3). If $X$ doesn't satisfy the magic rule we are in situation (1). Thus it remains to prove the following statement.

Let $(X,d)$ be metric space having $n=4,5,6$ or $7$ points, having at least two colors and satisfying the magic rule. Then (1), (2) or (4) holds.

For $n=4$ the only spaces which are left are the square and the rectangle. The duplex rule applies to the square and the star rule applies to the rectangle.

For $n=5,7$ lemma 4.1 shows that all color components must be $n$-cycles. The cycle rule applies and shows that we are in situation (2).

For $n=6$ lemma 4.1 shows that all color components must be among the following 6 graphs: 6-cycle, 3 segments, 2 triangles, and their complements. If there is a 6-cycle the cycle rule applies and (2) must hold. The same happens whenever the complement of the 6-cycle is present. In the remaining cases, either the duplex rule applies and (4) must hold, or the star rule applies and we are in situation (2).
\end{proof}

In the $n=8$ case the situation is a bit more complicated. There are 14 homogenous graphs with 8 vertices, and an analogue of lemma 4.1 holds. Techniques in this paper apply to most of them, and there are 3 graphs left. First is the cube, which can be shown to correspond to a tensor product between 2 Temperley-Lieb algebras. Second is a graph which looks like a wheel with 8 spokes, which gives $\c (D_8)$. Third is the graph formed by 2 squares. This corresponds to a Fuss-Catalan algebra on 3 colors, but the proof is quite long.

As for the passage from graphs to metric spaces, or colored graphs, this is done by an easy application of the decoloration rule, as in proof of theorem 4.1. However, there is a problem at $n=8$, with one space left, which is formed by 2 rectangles. We don't know yet how to compute invariants of the associated Hopf algebra.

For $n=9$ an analogue of lemma 4.1 holds one again, but the proof is long. All graphs give $\c (D_9)$ or correspond to Fuss-Catalan algebras, except for an exceptional graph, which looks like a torus. So far, we have no results about this graph.

\end{document}